\documentclass{article}

\usepackage{amsmath}
\usepackage{amsfonts}
\usepackage{IEEEtrantools}
\usepackage{stmaryrd}
\usepackage{tcolorbox}
\usepackage{hyphenat}

\usepackage{subcaption}

\newcommand{\mvecfun}[1]{ \underline{#1}}
\newcommand{\mtenfun}[1]{ \underline{\underline{#1}}}
\newcommand{\lbracket}[0]{ {[\hspace{0.02cm}} }
\newcommand{\rbracket}[0]{ {\hspace{0.02cm]}} }

\newcommand*{\ldblbrace}{\{\mskip-5mu\{}
\newcommand*{\rdblbrace}{\}\mskip-5mu\}}

\newcommand{\jump}[1]{ \lbracket {#1} \rbracket}
\newcommand{\jjump}[1]{ \llbracket {#1} \rrbracket}
\newcommand{\mean}[1]{ \{ {#1} \}}
\newcommand{\mmean}[1]{ \ldblbrace {#1} \rdblbrace}

\usepackage{authblk}

\newcommand{\diff}[0]{ \,\mathrm{d}}

\usepackage{setspace}
\date{}
\author[1,2]{Marco Favino}

\affil[1]{Faculty of Mathematics and Computer Science, UniDistance Suisse}
\affil[2]{Euler Institute, Universit\`a della Svizzera italiana}

\title{Fractures and thin heterogeneities\\as Robin-Wentzell interface conditions}

\begin{document}

\maketitle

\section*{Abstract}

We formally derive interface conditions for modeling fractures in Darcy flow problems and, more generally, thin inclusions in heterogeneous diffusion problems expressed as the divergence of a flux. Through a formal integration of the governing equations within the inclusions, we establish that the resulting interface conditions are of Wentzell type for the flux jump and Robin type for the flux average. Notably, the flux jump condition is unconventional, involving a tangential diffusion operator applied to the average of the solution across the interface.

The corresponding weak formulation is introduced, offering a framework that is readily applicable to finite element discretizations. Extensive numerical validation highlights the robustness and versatility of the proposed modeling technique. The results demonstrate its effectiveness in accommodating a wide range of material properties, managing networks of inclusions, and naturally handling fractures with varying apertures -- all without requiring an explicit geometric representation of the fractures.

\section{Introduction}

A wide range of applications throughout the natural and engineering sciences relies on computer simulations involving the numerical solution of heterogeneous differential problems featuring thin inclusions.  
By \textit{thin}, we mean that one dimension of the inclusion is orders-of-magnitude smaller than the characteristic size of the embedding background, while its material properties may also differ significantly from those of the background itself.

In the Earth, environmental, and engineering sciences,
a key computational challenge involving thin inclusions is the simulation of fluid flow in fractured porous media~\cite{martin2005modeling,flemisch2018benchmarks,berre2019flow,BenchPaper,favino2021}.  
Fractures are thin inclusions embedded within the rock mass, with their smallest dimension commonly referred to as the \emph{aperture}.  
The fluid flow problem is governed by two equations:  
the mass balance equation, involving the divergence of the volumetric flux,  
and Darcy's law, which relates the volumetric flux to the gradient of the fluid pressure via the fluid mobility.  
Despite their thin apertures, fractures play a dominant role in fluid flow dynamics.  
Conductive fractures, with a fluid mobility much larger than that of the embedding background, act as conduits,  
while blocking fractures, with a much smaller fluid mobility, act as barriers.

Mathematical models that explicitly account for the small
but finite aperture of fractures are referred to as equi-dimensional.
However, these models face significant challenges, particularly with regard to mesh generation.  
For realistic networks containing hundreds of fractures,
automated mesh generation is often infeasible
and instead requires extensive human intervention.
Even when successful, the thus generated meshes often include highly elongated elements,
which compromise the accuracy of the numerical solutions.  
Approximate strategies based on adaptive mesh refinement
have shown promise in two-dimensional settings~\cite{favino2020fully,favino2021},
but their computational cost becomes prohibitive for three-dimensional problems.

To address these challenges, reduced modeling techniques such as hybrid-dimensional representations of fractures are often employed~\cite{martin2005modeling,flemisch2018benchmarks,berre2019flow,BenchPaper}.
In hybrid\hyp{}dimensional models, fractures are represented as objects of a lower geometrical dimension than the embedding background.  
These techniques treat the background and fractures as distinct domains, with two sets of equations defined in each, coupled through appropriate continuity conditions.  
Treating fractures and the background as distinct domains allows for a significant simplification of mesh generation, albeit at the cost of introducing a small modeling error due to the hybrid-dimensional representation.

Two primary families of hybrid-dimensional models are recognized in the literature.
The first family derives from averaging the governing equations across the fracture aperture~\cite{martin2005modeling,BenchPaper},
leading to reduced-dimensional representations of the fractures.
The second family is based on a theoretical framework akin to that of immersed methods~\cite{PESKIN1972252,GRIFFITH200575},
resulting in the so-called \textit{continuous pressure models}~\cite{helmig1997multiphase,ZULIAN2022110773}.

Two main families of hybrid-dimensional models exist.  
The first is based on averaging the equations across the fracture aperture~\cite{martin2005modeling,BenchPaper}, while the second relies on a theoretical framework akin to that of immersed methods~\cite{PESKIN1972252,GRIFFITH200575},
resulting in the so-called \textit{continuous pressure models}~\cite{helmig1997multiphase,ZULIAN2022110773}.
Averaging-based models are in general effective for both conductive and blocking fractures but require defining equations at fracture intersections and intersections of fracture intersections using a recursive approach.  
Continuous pressure models naturally handle fracture intersections but are primarily suited for conductive fractures.  
Both approaches have been widely adopted for simulating fluid flow in fractured media.
For a detailed comparison of various methods employing hybrid-dimensional representations of fractures,
see, for example,~\cite{flemisch2018benchmarks,BenchPaper} and the references therein.

Thin inclusions are also present in relevant applications in other fields.  
For example, one problem in biomechanics that presents a mathematical structure similar to fluid flow in fractured porous media is cell migration through thin membranes~\cite{chaplain2019derivation}. In these processes, cells navigate barriers, such as other cells, cell-cell junctions, and extracellular matrices with varying density and composition.  
For corresponding numerical simulations, in such applications, alternative reduced modeling techniques for thin inclusions based on interface problems are typically preferred.

\textit{Interface problems} are characterized by solutions that may exhibit discontinuities across an interface, where appropriate conditions are imposed
to govern the jump and the average of either the flux or the solution~\cite{gong2008immersed,mazzucato2013nonconforming,burman2015cutfem,kwak2017new,le2017elliptic}. In the context of reduced modeling techniques, these interface conditions naturally emerge from the formal integration of the equations within the heterogeneity across its aperture~\cite{kedem1963permeability,chaplain2019derivation}, in a manner analogous to the averaging-based hybrid-dimensional models. Replacing the inclusion with interface conditions is particularly advantageous, as it allows to consider only the equations in the embedding background without the need to explicitly model the inclusions.

In this study, building upon the approach presented in~\cite{chaplain2019derivation},  
we formally derive conditions for modeling fractures as interfaces in Darcy flow problems.  
However, in contrast to their findings,
where the flux jump was found to be zero,
we demonstrate that the flux jump corresponds to a tangential diffusion operator
applied to the average pressure across the interface.
This condition can be interpreted as a generalization of Wentzell-type boundary conditions~\cite{venttsel1959boundary,bennequin2016optimized} to interfaces.
The second condition, on the other hand, is of Robin type for flux averages.  
 Incidentally, in formalizing this derivation, we establish a comprehensive framework for the possible interface conditions  
that can be imposed in interface problems - a framework that, to the best of our knowledge, is currently unavailable in the literature.
We introduce the corresponding weak formulation, which is readily applicable for finite element discretizations.  
The formulation is both simple and elegant, incorporating, in addition to the standard domain integrals, two types of interface integrals related to the flux jump and average conditions.  
The proposed model is validated through four numerical test cases involving fractured media, demonstrating its applicability to both conductive and blocking fractures, as well as its ability to naturally handle fracture intersections with varying apertures -- all without requiring an explicit fracture mesh.  
 
Although this work focuses on Darcy flow in fractured media, 
the proposed methodology is applicable to any diffusion problem involving thin heterogeneities 
and can be naturally extended to more general problems formulated in divergence form.

The manuscript is organized as follows. In Section~\ref{sec:1d}, we introduce the modeling strategy using a one-dimensional example. Section~3 provides the necessary mathematical preliminaries, including the geometric framework, the definitions of tangent and normal operators at the interface, and the formulation of Wentzell boundary conditions. In Section~4, we present a general framework for possible interface conditions in interface problems. Section~5 formally demonstrates how the interface problem can approximate the original problem for the case of a single thin inclusion. Section~6 extends the results to the case of a network of inclusions and introduces the corresponding weak formulation. Section~7 presents numerical validation of the proposed model through a series of test cases. Finally, Section~8 concludes the manuscript and discusses potential future directions.

\section{A One-Dimensional Example}
\label{sec:1d}

\subsection{Geometry of the problem}

\begin{figure}
\begin{subfigure}{.5\textwidth}
\centering
\includegraphics[scale=0.3]{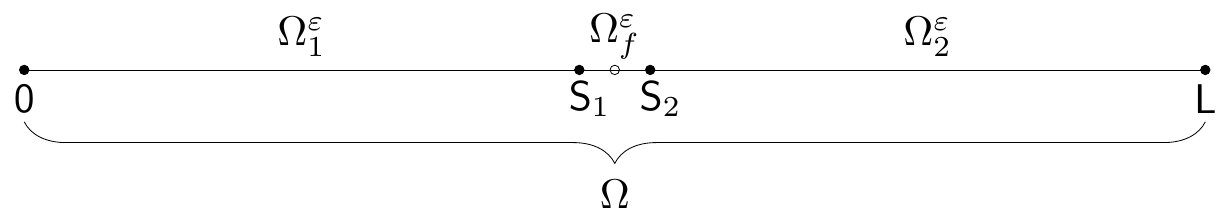}
\caption{Decomposition of \(\Omega\) for the equi-dimensional heterogeneous problem.}
\label{fig:decomposition1D1}
\end{subfigure}%

\begin{subfigure}{.5\textwidth}
\centering
\includegraphics[scale=0.3]{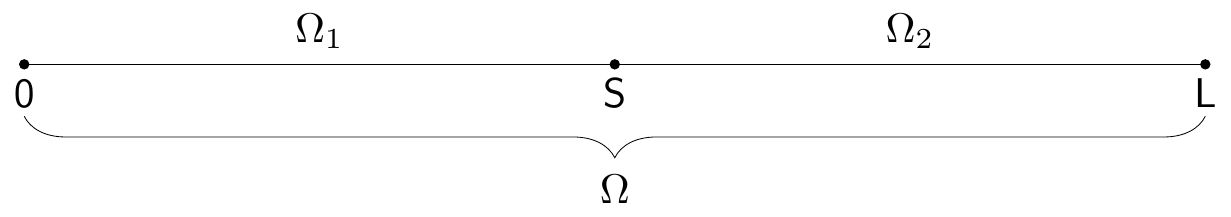}
\caption{Decomposition of \(\Omega\) for the interface problem.}
\label{fig:decomposition1D2}
\end{subfigure}
\caption{Illustration of domain decompositions for the heterogeneous and interface problems.}
\label{fig:decomposition1D}
\end{figure}

Let \(\Omega = (\mathsf{0}, \mathsf{L})\) be an interval in \(\mathbb{R}\)
and \(\mathsf{S}\) be a point in \(\Omega\).
For a positive constant \(\varepsilon \ll \mathsf{L}\) (referred to as the aperture),  
we define \(\mathsf{S}_1 = \mathsf{S} - \varepsilon/2\) and \(\mathsf{S}_2 = \mathsf{S} + \varepsilon/2\).  
The subdomain \(\Omega_f^\varepsilon = (\mathsf{S}_1, \mathsf{S}_2)\) represents a small inclusion in \(\Omega\).
The presence of \(\Omega_f^\varepsilon\) naturally introduces the following decomposition:
\begin{equation}
\label{eq:equidec}
\overline{\Omega} = \overline{\Omega}_1^\varepsilon \cup \overline{\Omega}_f^\varepsilon \cup \overline{\Omega}_2^\varepsilon,
\end{equation}
where:
\[
\Omega_1^\varepsilon = (\mathsf{0}, \mathsf{S}_1), \quad 
\Omega_2^\varepsilon = (\mathsf{S}_2, \mathsf{L}).
\]
This decomposition explicitly represents the thin inclusion as a subdomain of \(\Omega\)
(Figure~\ref{fig:decomposition1D1})
and will be referred to as equi-dimensional.

In addition, we introduce the decomposition:
\begin{equation}
\label{eq:hybriddec}
\overline{\Omega} = \overline{\Omega}_1 \cup \overline{\Omega}_2,
\end{equation}
where:
\[
\Omega_1 = (\mathsf{0}, \mathsf{S}), \quad \Omega_2 = (\mathsf{S}, \mathsf{L}).
\]
This decomposition will be used when the thin inclusion is replaced by an interface.
This interface-problem decomposition is illustrated in Figure~\ref{fig:decomposition1D2}.

Given a function \(v : \Omega \to \mathbb{R}\), we denote its restriction to \(\Omega_i^\varepsilon\)
by \(v_i^\varepsilon\), with \(i \in \{ 1,2,f \}\), and its restriction to \(\Omega_i\) by \(v_i\), with \(i \in \{ 1,2 \}\). 
The superscript \(\varepsilon\) is used to explicitly denote the dependence on the small inclusion, distinguishing it from the domain where interface conditions are applied.

\subsection{Heterogeneous Problems with Thin Inclusions}

We are interested in solving:
\begin{align}
    u' &= 0, \label{eq:balance1d} \\
    u &= -k p', \label{eq:constitutive1d}
\end{align}
with suitable boundary conditions, where \((\,\cdot\,)^\prime\) denotes the derivative.  
The parameter \(k\) is a diffusion coefficient that takes the value \(k_i\) in \(\Omega_i^\varepsilon\) for \(i \in \{1, 2, f\}\).  

In the context of fluid flow in fractured media:
\begin{itemize}
    \item Equation~\eqref{eq:balance1d} describes the conservation of mass in a porous medium, where \(u\) denotes the volumetric flux.
    \item Equation~\eqref{eq:constitutive1d} is Darcy's law, which relates the derivative of the pressure \(p\) with the flux \(u\), and \(k\) represents the fluid mobility, defined as the ratio between the permeability of the porous medium and the fluid viscosity.
\end{itemize}

In a more general setting,  
Equation~\eqref{eq:balance1d} represents a balance law with \(u\) as the associated flux,  
while Equation~\eqref{eq:constitutive1d} is the closure relation connecting \(u\) with a primary variable \(p\).

Considering the decomposition~\eqref{eq:equidec}, the differential problem~\eqref{eq:balance1d} and~\eqref{eq:constitutive1d} can be written as:
\begin{alignat}{2}
(u_i^\varepsilon)^\prime & = 0  & \quad \text{in } \Omega _{i}^\varepsilon, & \quad  i\in\{ 1,2,f \},
\label{eq:balance1di}\\
u_i^\varepsilon & = -k_i(p_i^\varepsilon)^\prime  & \quad \text{in } \Omega _{i}^\varepsilon, & \quad  i\in\{ 1,2,f \},
\label{eq:constitutive1di}\\
u_i^\varepsilon(\mathsf{S}_i) & = u_f^\varepsilon(\mathsf{S}_i) & & \quad i\in\{ 1,2 \},
\label{eq:continuityui}\\
p_i^\varepsilon(\mathsf{S}_i) & = p_f^\varepsilon(\mathsf{S}_i) & & \quad i\in\{ 1,2 \}.
\label{eq:continuitypi}
\end{alignat}
Equations~\eqref{eq:balance1di} and~\eqref{eq:constitutive1di} are the restrictions on \(\Omega_i^{\varepsilon}\)
of the balance law~\eqref{eq:balance1d} and the constitutive equation~\eqref{eq:constitutive1d}.
Equations~\eqref{eq:continuityui} and~\eqref{eq:continuitypi} represent the continuity conditions at \(\mathsf{S}_i\) for \(u\) and \(p\), respectively.

\subsection{Thin Inclusions as Interface Conditions}

Our goal is to derive a suitable interface problem where the equations in~$\Omega_f^\varepsilon$
are replaced with interface conditions, providing an approximation of the original heterogeneous problem.
We proceed as follows:

\smallskip

\noindent
First, we integrate the balance equation~\eqref{eq:balance1di} and the closure relation~\eqref{eq:constitutive1di} over $\Omega_f^\varepsilon$.
While the integrals of derived quantities can be directly evaluated,
the integral of $u_f^\varepsilon$ is approximated using a trapezoidal quadrature rule,
leading to:
\begin{alignat}{2}
u_f^\varepsilon(\mathsf{S}_2) - u_f^\varepsilon(\mathsf{S}_1) &=  0, \label{eq:step1u} \\[2mm]
\frac{\varepsilon}{2} \left( u_f^\varepsilon (\mathsf{S}_1) + u_f^\varepsilon (\mathsf{S}_2) \right)  &=  -k_f \left( p_f^\varepsilon(\mathsf{S}_2) - p_f^\varepsilon(\mathsf{S}_1) \right). \label{eq:step1p}
\end{alignat}
The use of a trapezoidal quadrature rule is equivalent to assuming that $u_f^\varepsilon$ is linear across $\Omega_f^\varepsilon$.

\smallskip

\noindent
Next, we apply the continuity conditions~\eqref{eq:continuityui} and~\eqref{eq:continuitypi} at $\mathsf{S}_i$:
\begin{alignat}{2}
u_2^\varepsilon(\mathsf{S}_2) - u_1^\varepsilon(\mathsf{S}_1) &=  0, \label{eq:step2u} \\[2mm]
\frac{\varepsilon}{2} \left( u_1^\varepsilon (\mathsf{S}_1) + u_2^\varepsilon (\mathsf{S}_2) \right)  &=  -k_f \left( p_2^\varepsilon(\mathsf{S}_2) - p_1^\varepsilon(\mathsf{S}_1) \right). \label{eq:step2p}
\end{alignat}

\smallskip

\noindent
Finally, we take the limits $\mathsf{S}_i \to \mathsf{S}$, yielding the interface conditions:
\begin{alignat}{2}
u_2(\mathsf{S}) - u_1(\mathsf{S}) &=  0, \label{eq:step3u} \\[2mm]
\frac{\varepsilon}{2} \left( u_1 (\mathsf{S}) + u_2 (\mathsf{S}) \right)  &=  -k_f \left( p_2(\mathsf{S}) - p_1(\mathsf{S}) \right), \label{eq:step3p}
\end{alignat}
where we replaced $u_i^\varepsilon$ and $p_i^\varepsilon$ with $u_i$ and $p_i$, respectively.

Hence, the heterogeneous problem composed of Equations~\eqref{eq:balance1di}, \eqref{eq:constitutive1di},
\eqref{eq:continuityui}, and \eqref{eq:continuitypi}
has been rewritten as the following interface problem:
\begin{alignat}{2}
u_i^\prime &= 0  & \quad \text{in } \Omega _{i}, & \quad  i \in \{ 1,2 \}, \label{eq:balance1di2} \\[1mm]
u_i &= -k_i p_i^\prime  & \quad \text{in } \Omega _{i}, & \quad  i \in \{ 1,2 \}, \label{eq:constitutive1di2} \\[1mm]
-\jump{u} &=  0, \label{eq:interfaceu} \\[1mm]
 -\mean{u} &=  \frac{k_f}{\varepsilon} \jump{p}, \label{eq:interfacep}
\end{alignat}
where Equations~\eqref{eq:interfaceu} and~\eqref{eq:interfacep} represent the interface conditions,
written using the standard notations for the jump and average of a function at $\mathsf{S}$, defined respectively as:
\[
\jump{v} = v_2(\mathsf{S}) - v_1(\mathsf{S}), \quad 
\mean{v} = \frac{1}{2} \left( v_1(\mathsf{S}) + v_2(\mathsf{S}) \right).
\]

\subsection{Variational and Finite Element Formulation}

We now consider, as an example, the interface problem defined by Equations~\eqref{eq:balance1di2}- \eqref{eq:interfacep}, completed with the following Dirichlet and Neumann boundary conditions:
\[
u(\mathsf{0}) = u_1(\mathsf{0}) = h, \quad p(\mathsf{L}) = p_2(\mathsf{L}) = 0.
\]
We introduce the function space:
\[
V = \{ v \in L^2(\Omega) \, : \, v_i \in H^1(\Omega_i) \, \text{ for } i = 1, 2, \text{ and } v(\mathsf{L}) = 0 \}.
\]

The variational formulation is derived by multiplying the balance equation~\eqref{eq:balance1di2} by a test function \(v \in V\), integrating over the domain, and applying integration by parts. This yields:
\begin{multline}
\label{eq:weakform1d}
- \int_{\mathsf{0}}^{\mathsf{S}} u_1 \, v_1^\prime \, \mathrm{d}x 
- \int_{\mathsf{S}}^{\mathsf{L}} u_2 \, v_2^\prime \, \mathrm{d}x \\ 
- \jump{u} \mean{v} - \mean{u} \jump{v} 
+ u_2(\mathsf{L}) v_2(\mathsf{L}) - u_1(\mathsf{0}) v_1(\mathsf{0}) = 0.
\end{multline}

Replacing the constitutive law~\eqref{eq:constitutive1di2}, the interface conditions~\eqref{eq:interfaceu} and~\eqref{eq:interfacep}, the Neumann boundary condition, and noting that \(v_2(\mathsf{L}) = 0\), we obtain the weak formulation:

\smallskip

\noindent
{Find} \(p \in V\) {such that}
\begin{equation}
\label{eq:weak1D}
\int_{\mathsf{0}}^{\mathsf{S}} k_1 \, p_1^\prime \, q_1^\prime \, \mathrm{d}x 
+ \int_{\mathsf{S}}^{\mathsf{L}} k_2 \, p_2^\prime \, q_2^\prime \, \mathrm{d}x 
+ \frac{k_f}{\varepsilon} \jump{p} \jump{q} 
= h \, q_1(\mathsf{0}) \quad \forall q \in V.
\end{equation}

The weak formulation above can be directly used for finite element discretization. Proceeding in the standard way, we introduce two meshes for \(\Omega_1\) and \(\Omega_2\), define the usual finite-dimensional subspace \(V_h\), and state the problem~\eqref{eq:weak1D} in \(V_h\). Notably, at \(\mathsf{S}\) the degrees of freedom are duplicated, since the vertex \(\mathsf{S}\) exists in \(\Omega_1\) and \(\Omega_2\). This allows for a discontinuity in the solution at the interface.

\begin{figure}
\centering
\includegraphics[scale=0.4]{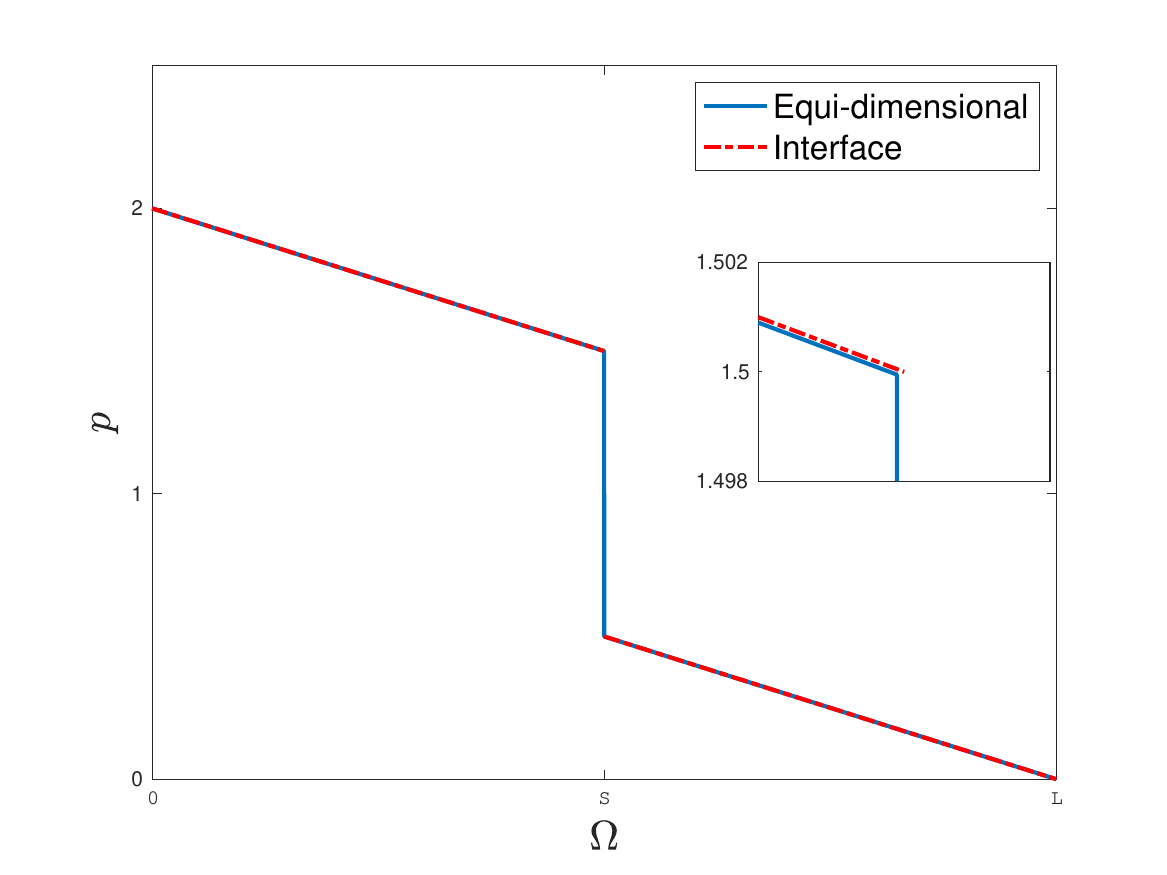}
\caption{Comparison between the solutions obtained with the heterogeneous and interface problems.}
\label{fig:comparison1D}
\end{figure}

In Figure~\ref{fig:comparison1D}, we compare the solutions obtained from the heterogeneous and the interface problems using a linear finite element method for \(\mathsf{L} = k_1 = k_2 = h = 1\), \(\varepsilon = k_f = 10^{-4}\), and \(\mathsf{S} = 0.5\). The results show good macroscopic agreement between the two models. 
We observe that \(p\) exhibits a significant drop within the small inclusion due to the small value of \(k_f\), and the interface problem successfully reproduces this behavior through the interface conditions. However, closer inspection in the zoomed-in view reveals small differences, indicating modeling errors introduced by the approximation.

\section{Mathematical Preliminaries}

\subsection{Decompositions of the domain}
\label{sec:dec}

\begin{figure}
\begin{subfigure}{1\textwidth}
\centering
\includegraphics[scale=2]{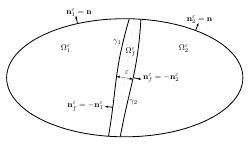}
\caption{Decomposition of \(\Omega\) for the equi-dimensional heterogeneous problem.}
\label{fig:decomposition2D1}
\end{subfigure}%

\begin{subfigure}{1\textwidth}
\centering
\includegraphics[scale=2]{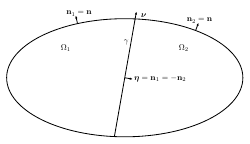}
\caption{Decomposition of \(\Omega\) for the interface problem.}
\label{fig:decomposition2D2}
\end{subfigure}
\caption{Illustration of the decompositions of $\Omega$ for the heterogeneous and interface problems.}
\label{fig:decomposition2D}
\end{figure}

Throughout this manuscript, we denote a bounded domain with Lipschitz boundary in \(\mathbb{R}^d\), with \(d \in \{2,3\}\), by \(\Omega\). 
The Cartesian coordinates of \(\Omega\) are denoted by \(x\) and \(y\) for \(d=2\), and \(x\), \(y\), and \(z\) for \(d=3\). 
Points in \(\Omega\) are represented using uppercase sans serif letters, such as \(\mathsf{A}\). 
We also write \(\mathsf{A} = (x, y)\) or \(\mathsf{A} = (x, y, z)\) to explicitly specify their coordinates.

The boundary of \(\Omega\) is denoted by \(\Gamma = \partial \Omega\), 
and the outward unit normal to \(\Gamma\) is represented by \(\mathbf{n}\). 
We assume that \(\Gamma\) is decomposed into two subsets, \(\Gamma^{p}\) and \(\Gamma^{\mathbf{u}}\), such that:
\[
\overline{\Gamma}^p \cup \overline{\Gamma}^{\mathbf{u}} = \Gamma, \quad \text{and} \quad \mathring{\Gamma}^p \cap \mathring{\Gamma}^{\mathbf{u}} = \emptyset,
\]
where \(\Gamma^p\) is assumed to have a positive measure, i.e., \(|\Gamma^p| > 0\).

Let \(\gamma \subset \Omega\) be a smooth interface, i.e., a curve for \(d=2\) or a surface for \(d=3\),  
which separates \(\Omega\) into two connected subdomains.  
We denote a unit normal to \(\gamma\) by \(\boldsymbol{\eta}\) and the outward-pointing unit normal vector to \(\partial \gamma\) by \(\boldsymbol{\nu}\).  

We introduce the subdomain \(\Omega_f^\varepsilon \), representing a thin inclusion in $\Omega$, defined as:
\[
\Omega_f^\varepsilon = \left\{ \mathsf{X} \in \Omega \, : \, \mathsf{X} = \mathsf{S} + r \frac{\varepsilon(\mathsf{S})}{2} \boldsymbol{\eta}, \; \mathsf{S} \in \gamma, \; -1 < r < 1 \right\},
\]
where \(\varepsilon(\mathsf{S})\) is a sufficiently smooth function mapping \(\gamma\) to \(\mathbb{R}^+\), and represents the aperture of the inclusion at \(\mathsf{S} \in \gamma\) along the normal direction \(\boldsymbol{\eta}\).

The subset \(\Omega_f^\varepsilon\) naturally leads to the \emph{equi\hyp{}dimensional} decomposition of \(\Omega\)
into three subdomains, as depicted in Figure~\ref{fig:decomposition2D1}, such that:
\begin{equation}
\label{eq:equiDec}
\overline{\Omega}= \overline{\Omega}_1^\varepsilon  \cup \overline{\Omega}_f^\varepsilon \cup \overline{\Omega}_2^\varepsilon.
\end{equation}
As a result, the boundary \(\partial \Omega\) can be divided into the subsets \(\Gamma_i^\varepsilon = \Gamma \cap \overline{\Omega}_i^\varepsilon\), for \(i \in \{1, 2, f\}\). We also identify the sets \(\gamma_i = \overline{\Omega}_i^\varepsilon \cap \overline{\Omega}_f^\varepsilon\), with \(i \in \{1, 2\}\). Denoting the outward normal to \(\Omega_i^\varepsilon\) by \(\mathbf{n}_i^\varepsilon\), it follows that \(\mathbf{n} = \mathbf{n}_i^\varepsilon\) on \(\Gamma_i^\varepsilon\), and \(\mathbf{n}_f^\varepsilon = -\mathbf{n}_i^\varepsilon\) on \(\gamma_i\).

Instead, the interface $\gamma$ naturally leads to the \emph{interface problem} decomposition of \(\Omega\)
into two subdomains, as depicted in Figure~\ref{fig:decomposition2D2}, such that:
\begin{equation}
\label{eq:interDec}
\overline{\Omega} = \overline{\Omega}_1 \cup \overline{\Omega}_2 \quad \text{and} \quad \mathring{\Omega}_1 \cap \mathring{\Omega}_2 = \emptyset,
\end{equation}
where \(\Omega_i\) is such that \(\Omega_i \cap \Omega_i^\varepsilon \neq \emptyset\).
As a result, the boundary \(\partial \Omega\) can be divided into the subsets \(\Gamma_i = \Gamma \cap \overline{\Omega}_i\), for \(i \in \{1, 2\}\) and \(\Gamma_i^{\beta} = \Gamma^\beta \cap \overline{\Omega}_i\),  
with \(i \in \{1,2\}\) and \(\beta \in \{p, \mathbf{u}\}\).
Denoting the outward normal to \(\Omega_i\) by \(\mathbf{n}_i\), it follows that \(\mathbf{n} = \mathbf{n}_i\) on \(\Gamma_i\), and \(\mathbf{n}_1 = -\mathbf{n}_2\) on \(\gamma\).  
Without loss of generality, we assume that \(\boldsymbol{\eta} = \mathbf{n}_1\) on \(\gamma\).

For a sufficiently smooth scalar function \(q\) and vector function \(\mathbf{v}\) defined over \(\Omega\),  
we denote their restriction to \(\Omega_i\) by \(q_i = q |_{\Omega_i}\) and \(\mathbf{v}_i = \mathbf{v} |_{\Omega_i}\), respectively, and their restrictions to \(\Omega_i^\varepsilon\) by \(q_i^\varepsilon = q |_{\Omega_i^\varepsilon}\) and \(\mathbf{v}_i^\varepsilon = \mathbf{v} |_{\Omega_i^\varepsilon}\), respectively.
We will make use of the standard notation used to denote the jump and average of scalar and vector functions at the interface \(\gamma\), i.e.:
\[
\begin{aligned}
\jump{q} & = q_2{}_{| \gamma} - q_1{}_{| \gamma}, & \quad & \jjump{\mathbf{v}} = (\mathbf{v}_2{}_{| \gamma} - \mathbf{v}_1{}_{| \gamma}) \cdot \boldsymbol{\eta}, \\
\mean{q} & = \frac{q_1{}_{| \gamma} + q_2{}_{| \gamma}}{2}, & \quad & \mmean{\mathbf{v}} = \frac{(\mathbf{v}_1{}_{| \gamma} + \mathbf{v}_2{}_{| \gamma}) \cdot \boldsymbol{\eta}}{2}.
\end{aligned}
\]



Given a set \(\Xi\), we adopt the usual notations \(L^s(\Xi)\) and \(H^s(\Xi)\) for the Lebesgue and Sobolev spaces of order \(s\), respectively.

\subsection{Tangential and Normal Differential Operators}

Throughout this manuscript, we denote the gradient and divergence operators by \(\nabla\) and \(\nabla \cdot\), respectively.  
For a scalar function \(q\) and a vector function \(\mathbf{v}\), they are written as \(\nabla q\) and \(\nabla \cdot \mathbf{v}\).  

Let \(\Sigma\) represent either the interface \(\gamma\) or a subset of \(\Gamma\),  
and let \(\boldsymbol{\eta}\) denote its unit normal vector.  
We define the projection operator onto the plane tangent to \(\Sigma\) as  
\[
\mathbf{P} = \mathbf{I} - \boldsymbol{\eta} \otimes \boldsymbol{\eta},
\]  
where \(\mathbf{I}\) is the identity tensor.  

Using this projection, the tangential component of a vector field \(\mathbf{v}\) is given by  
\[
\mathbf{v}_\tau = \mathbf{P} \mathbf{v},
\]  
and its normal component by  
\[
v_\eta = \mathbf{v} \cdot \boldsymbol{\eta}.
\]  
Thus, the projection of \(\mathbf{v}\) along \(\boldsymbol{\eta}\) is  
\[
\mathbf{v}_\eta = v_\eta \boldsymbol{\eta}.
\]  
Consequently, the decomposition of \(\mathbf{v}\) into its tangential and normal components is given by
\[
\mathbf{v} = \mathbf{v}_\tau + \mathbf{v}_\eta.
\]

We now define the tangential and normal gradient operators, denoted by \(\nabla_\tau\) and \(\nabla_\eta\),  
as well as the tangential and normal divergence operators, denoted by \(\nabla_\tau \cdot\) and \(\nabla_\eta \cdot\).  
For gradients, we have the following identities:  
\[
\nabla_\tau q = \mathbf{P} \nabla q \quad \text{and} \quad \nabla_\eta q = \frac{\partial q}{\partial \eta} \boldsymbol{\eta},
\]  
where \(\frac{\partial q}{\partial \eta}\) denotes the directional derivative along \(\boldsymbol{\eta}\).  
Thus, the gradient can be decomposed as  
\[
\nabla q = \nabla_\tau q + \nabla_\eta q.
\]

For divergences, the decomposition is given by:  
\begin{equation}
\label{eq:divEq}
\nabla \cdot \mathbf{v} = \nabla_\tau \cdot \mathbf{v} + \nabla_\eta \cdot \mathbf{v} = \nabla_\tau \cdot \mathbf{v}_\tau + \nabla_\eta \cdot \mathbf{v}_\eta.
\end{equation}
Finally, we note that  
\begin{equation}
\label{eq:normalDiv}
\nabla_\eta \cdot \mathbf{v}_\eta = \frac{\partial v_\eta}{\partial \eta}.
\end{equation}

\subsection{Model Problem and Wentzell Boundary Conditions}

The fluid flow in porous media is described by the following system of equations:
\begin{alignat}{2}
-\nabla \cdot \mathbf{u} &= 0 & \quad & \text{in } \Omega, \label{eq:balance} \\
\mathbf{u} &= -k \nabla p & \quad & \text{in } \Omega, \label{eq:closure}
\end{alignat}
where Equation~\eqref{eq:balance} represents the conservation of mass for an incompressible fluid  
in a porous medium, and Equation~\eqref{eq:closure}, known as Darcy's law,  
arises from the balance of linear momentum.  

Here, \(\mathbf{u}\) is the volumetric flux (physically a velocity),  
\(p\) is the fluid pressure,  
and \(k\) is the fluid mobility, defined as the ratio of the permeability of the porous medium to the fluid viscosity.  
As we are interested in heterogeneous problems, we assume that \(k = k(X)\) is spatially varying in \(\Omega\),  
and bounded and strictly positive. Specifically, there exist constants \(k_{\min}\) and \(k_{\max}\) such that:  
\[
k_{\min} \leq k(X) \leq k_{\max}, \qquad \forall X \in \Omega.
\]  

The system described by Equations~\eqref{eq:balance} and~\eqref{eq:closure} is a  
commonly encountered formulation of diffusion problems. 
In this context, Equation~\eqref{eq:balance} enforces a balance principle within the domain,  
where \(\mathbf{u}\) represents the associated flux.  
Meanwhile, Equation~\eqref{eq:closure} serves as a closure relation, linking the flux \(\mathbf{u}\)  
to the primary variable \(p\).  
Hence, in the rest of the manuscript, we will generally refer to the conservation of mass  
as the balance equation,  
and to Darcy's law as the closure relation.

The so-called standard primal formulation is obtained by substituting the closure relation  
into the balance equation, resulting in the standard diffusion equation:  
\[
-\nabla \cdot (k \nabla p) = 0 \quad \text{in } \Omega.
\]

Problem~\eqref{eq:balance}-\eqref{eq:closure} must be supplemented with boundary conditions  
to ensure the existence and uniqueness of the solution.  
The type of boundary conditions that can be imposed on subsets of \(\Gamma\) becomes evident during the derivation of the corresponding weak formulation. Consider a test function \(q \in V \subset H^1(\Omega)\), where \(V\) will be specified later. Multiplying the balance equation~\eqref{eq:balance} by \(q\) and integrating over the domain, we obtain:

\[
\int_{\Omega} \nabla \cdot \mathbf{u} \, q \, \text{d}X = -\int_{\Omega} \mathbf{u} \cdot \nabla q \, \text{d}X + \underbrace{\int_{\Gamma} \mathbf{u} \cdot \mathbf{n} \, q \, \text{d}\sigma}_{\text{boundary integral}}.
\]
The boundary integral reveals the types of boundary conditions that can be imposed. Specifically, we can prescribe either the value of the pressure \(p\), which usually constrains the test function \(q\), or the value of the inward normal flux \(-\mathbf{u} \cdot \mathbf{n}\). When the pressure is prescribed, we have the {\bf Dirichlet} boundary condition:

\[
p = g \quad \text{on } \Gamma^{p}, \tag{1}
\]

In contrast, for the normal flux, the more general boundary condition commonly encountered in applications is:

\[
-\mathbf{u} \cdot \mathbf{n} = h - r p + \nabla_{\tau} \cdot (\kappa \nabla_\tau p) \quad \text{on } \Gamma^{\mathbf{u}}, \tag{2}
\]

where \(h\), \(r\geq0\), and \(\kappa\geq0\) are functions defined on \(\Gamma^{\mathbf{u}}\). The boundary condition on the inward normal flux can be categorized as follows:

\begin{itemize}
    \item \textbf{Wentzell} if \(\kappa > 0\), see, e.g.,~\cite{venttsel1959boundary,favini2002heat,kashiwabara2015well,bennequin2016optimized,LANCIA2017265}
    \item \textbf{Robin} if \(\kappa = 0\) but \(r > 0\),
    \item \textbf{Neumann} if \(\kappa = 0\) and \(r = 0\).
\end{itemize}

Wentzell conditions are quite non-standard as they involve a diffusion operator on \(\Gamma^{\mathbf{u}}\)  
and also affect the definition of the function space $V$.  
While in the case of Robin and Neumann conditions, the function space is given by
\[
V = \{ q \in H^1(\Omega) \, : \, q|_{\Gamma^p} = 0 \},
\]
in the case of Wentzell boundary conditions, we require
\[
V = \{ q \in H^1(\Omega) \, : \, q|_{\Gamma^{\mathbf{u}}} \in H^1(\Gamma^{\mathbf{u}}) \text{ and } q|_{\Gamma^p} = 0 \}.
\]

\section{A General Framework for Interface Conditions in Interface Problems}

Interface problems are defined on the decomposition~\eqref{eq:interDec} of $\Omega$,
where both \(p\) and the flux \(\mathbf{u}\)  
may exhibit discontinuities at the interface \(\gamma\). For Dirichlet and Neumann boundary conditions, the governing equations are:  
\begin{alignat}{2}
-\nabla \cdot \mathbf{u}_i &= 0 & \quad & \text{in } \Omega_i, \quad \text{for } i=1,2, \label{eq:balanceInt} \\
\mathbf{u}_i &= -k_i \nabla p_i & \quad & \text{in } \Omega_i, \quad \text{for } i=1,2, \label{eq:closureInt}\\
-\mathbf{u} \cdot \mathbf{n} &= h & \quad & \text{on } \Gamma^{\mathbf{u}}, \label{eq:NeuInt}\\
p &= g & \quad & \text{on } \Gamma^p. \label{eq:DirInt}
\end{alignat}

In addition to boundary conditions, interface conditions are essential to ensure the existence and uniqueness of the solution.  
In the literature, interface conditions are typically imposed on the jump of the solution, \(\jump{p}\),  
and the jump of the flux, \(\jjump{\mathbf{u}}\); see, for example,~\cite{gong2008immersed,mazzucato2013nonconforming,burman2015cutfem,kwak2017new,le2017elliptic}.  
However, much like boundary conditions, the possible interface conditions can be systematically derived  
from the interface integrals that emerge during the formulation of the associated weak problem.

To formalize this, we introduce the function space:  
\begin{equation}
\label{eq:defV}
\begin{array}{ l l  l }
V = \{ q \in L^2(\Omega) \, : & \text{\phantom{ii}i) } q_i \in H^1(\Omega_i), &\\
& \text{\phantom{i}ii) } {q}_{| \Gamma^p}=0, &\\
& \text{iii) } {q_i}_{| \Gamma_i^{\mathbf{u}}} \in H^1(\gamma), & i=1,2 \} \\
\end{array}
\end{equation}
and the function set \(U\), defined as \(V\), but where hypothesis ii) is replaced by \(q_{| \Gamma^p}=g\).  

By multiplying the balance law~\eqref{eq:balanceInt} by a test function \(q \in V\), integrating over \(\Omega\), and applying integration by parts, the left-hand side becomes:  
\begin{equation} \label{eq1}
\begin{split}
\sum_{i=1}^2 \int_{\Omega_i} \nabla \cdot \mathbf{u}_i \, q_i \diff \Omega & = \sum_{i=1}^2 \left( - \int_{\Omega_i} \mathbf{u}_i \cdot \nabla q_i \diff \Omega + \int_{\Gamma_i} \mathbf{u}_i \cdot \mathbf{n}_i \, q_i \diff \Gamma \right) \\
& \quad - \underbrace{\int_\gamma \jjump{\mathbf{u}} \, \mean{q} \diff \Gamma}_{\text{Flux Jump Integral}}  
- \underbrace{\int_\gamma \mmean{\mathbf{u}} \, \jump{q} \diff \Gamma}_{\text{Flux Average Integral}}.
\end{split}
\end{equation}  

Two interface integrals emerge:  
\begin{enumerate}
\item \textbf{Flux Jump Integral}, involving the product of the flux jump \(\jjump{\mathbf{u}}\) and the average test function \(\mean{q}\).  
\item \textbf{Flux Average Integral}, involving the product of the flux average \(\mmean{\mathbf{u}}\) and the test function jump \(\jump{q}\).  
\end{enumerate}

For these terms, possible interface conditions are:  
\begin{itemize}
\item For the \textbf{Flux Jump Integral}:  
\begin{align}
\mean{p} &= g^a, \quad \text{on } \gamma, \label{eq:dirmean}
\end{align}
or  
\begin{align}
-\jjump{\mathbf{u}} &= -\nabla_{\Gamma} \cdot \left(\kappa^j \nabla \mean{p} \right) + r^j \mean{p} - h^j, \quad \text{on } \gamma. \label{eq:wentjump}
\end{align}
\item For the \textbf{Flux Average Integral}:  
\begin{align}
\jump{p} &= g^j, \quad \text{on } \gamma, \label{eq:dirjump}
\end{align}
or  
\begin{align}
-\mmean{\mathbf{u}} &= -\nabla_{\Gamma} \cdot \left(\kappa^a \nabla \jump{p} \right) + r^a \jump{p} + h^a, \quad \text{on } \gamma. \label{eq:wentmean}
\end{align}
\end{itemize}

Equations~\eqref{eq:dirmean} and~\eqref{eq:dirjump} correspond to interface conditions of \textbf{Dirichlet} type,  
while Equations~\eqref{eq:wentjump} and~\eqref{eq:wentmean} are of type:
\begin{itemize}
\item \textbf{Wentzell}: If \(\kappa^\beta > 0\);  
\item \textbf{Robin}: If \(\kappa^\beta = 0\) and \(r^\beta > 0\);
\item \textbf{Neumann}: If \(\kappa^\beta = 0\) and \(r^\beta = 0\).  
\end{itemize}
Despite the generality of this formalism for interface conditions,  
we could not find any references addressing interface conditions of Wentzell type.  
The only work in which somewhat similar conditions are employed is~\cite{burman2019simple},  
where the conditions used were \(\jump{p} = 0\) and \(-\jjump{\mathbf{u}} = \nabla_\tau \cdot \nabla_\tau p\).

\subsection{Primal Weak Formulation of Interface Problems with Wentzell Interface Conditions}  

We consider an interface problem composed of Equations~\eqref{eq:balanceInt},~\eqref{eq:closureInt},~\eqref{eq:NeuInt}, and \eqref{eq:DirInt},  
completed by the interface conditions of Wentzell type for both the flux average and flux jump,  
i.e., Equations~\eqref{eq:wentjump} and~\eqref{eq:wentmean}.  
As we will see, this definition of an interface problem may be incomplete, as additional conditions may be needed depending on the geometry of \(\gamma\).  

Replacing the closure relation~\eqref{eq:closureInt}, the Neumann boundary condition~\eqref{eq:NeuInt},  
the interface conditions~\eqref{eq:wentjump} and~\eqref{eq:wentmean} into~\eqref{eq1}, 
and setting the resulting term equal to zero, we obtain: 
\begin{equation} \label{eq:weak_form}  
\begin{split}  
 \sum_{i=1}^2 \int_{\Omega_i} \nabla p_i \cdot \nabla q_i \diff \Omega - \int_{\Gamma^{\mathbf{u}}} h \, q \diff \Gamma & \\  
 - \int_\gamma \nabla_{\Gamma} \cdot \left(\kappa^j \nabla_{\Gamma} \mean{p} \right) \, \mean{q} \diff \Gamma  
- \int_\gamma \nabla_{\Gamma} \cdot \left(\kappa^a \nabla_{\Gamma} \jump{p} \right) \, \jump{q} \diff \Gamma & \\  
 + \int_\gamma r^j \mean{p} \, \mean{q} \diff \Gamma + \int_\gamma r^a \jump{p} \, \jump{q} \diff \Gamma & \\  
  - \int_\gamma h^j \, \mean{q} \diff \Gamma - \int_\gamma h^a \, \jump{q} \diff \Gamma & = 0. \\  
\end{split}  
\end{equation}  

In the second row of the equation above, diffusion terms for \(\jump{p}\) and \(\mean{p}\) along the interface \(\gamma\) appear.  
For these terms, we apply integration by parts over \(\gamma\).  
For instance, for the jump term, we have:  
\[  
- \int_\gamma \nabla_{\Gamma} \cdot \left(\kappa^a \nabla_{\Gamma} \jump{p} \right) \, \jump{q} \diff \Gamma =  
\int_\gamma \kappa^a \nabla_{\Gamma} \jump{p} \cdot \nabla_{\Gamma} \jump{q} \diff \Gamma  
- \int_{\partial \gamma} \kappa^a \nabla_{\Gamma} \jump{p} \cdot \boldsymbol{\nu} \, \jump{q} \diff l.  
\]  

The integral over \(\partial \gamma\) indicates the need for additional boundary conditions.  
However, over the portion of \(\partial \gamma\) intersecting \(\Gamma^p\), this term vanishes due to the test functions.  
For the parts intersecting \(\Gamma^{\mathbf{u}}\), we assume:  
\begin{equation}
\label{eq:hp}
\kappa^\beta \nabla_{\Gamma} \jump{p} \cdot \boldsymbol{\nu} = 0, \text{ with } \beta\in\{ a,j \}.
\end{equation}
This assumption is consistent with other hybrid-dimensional models for fractured media~\cite{martin2005modeling,boon2018robust,BenchPaper}.  

Hence, the primal weak formulation associated with the considered interface problem reads:  

~\\
\noindent  
Find \(p \in U\) such that  
\begin{equation} \label{eq:weak_form2}  
\begin{array}{r l c l }
\displaystyle \sum_{i=1}^2 & \displaystyle\int_{\Omega_i}  \nabla p_i \cdot \nabla q_i \diff \Omega & + & \\
&\displaystyle  \int_\gamma \kappa^m \nabla_{\Gamma} \jump{p} \cdot \nabla_{\Gamma} \jump{q} \diff \Gamma &+ &
\displaystyle \int_\gamma \kappa^a \nabla_{\Gamma} \mean{p} \cdot \nabla_{\Gamma} \mean{q} \diff \Gamma  \\
& \displaystyle \int_\gamma r^j \mean{p} \, \mean{q} \diff \Gamma &+& \displaystyle\int_\gamma r^a \jump{p} \, \jump{q} \diff \Gamma\\
= &  \displaystyle \int_{\Gamma^{\mathbf{u}}} h \, q \diff \Gamma & &\\
+ &  \displaystyle \int_\gamma h^j \, \mean{q} \diff \Gamma & +  & \displaystyle \int_\gamma h^a \, \jump{q} \diff \Gamma, \quad \forall q \in V.\\
\end{array}
\end{equation}  

Although specific results regarding the uniqueness and stability of the solution to this problem could not be found,  
a comparison with similar interface problems suggests that \(\kappa^\beta\) and \(r^\beta\) for \(\beta \in \{a, j\}\) can also be null  
if both \(\Gamma_i^p\) with \(i \in \{1, 2\}\) have positive measure.  
In other cases, at least one of them is expected to be nonzero.  
Finally, if both \(\kappa^\beta\) for \(\beta \in \{a, j\}\) are identically zero,  
the hypothesis iii) in the definition~\eqref{eq:defV} of \(V\) is no longer required,  
as there are no diffusion terms along the interface \(\gamma\).

\section{Thin inclusions as interface conditions}
\label{sec:derivation}

In this section, we formally derive a reduced modeling technique to address thin heterogeneities 
in fluid flow through fractured porous media.
By replacing thin heterogeneities with suitable interface conditions, 
we simplify the computational effort involved in mesh generation while preserving the essential physical characteristics of the flow.  

We aim to solve the problem defined by~\eqref{eq:balance}--\eqref{eq:closure}, defined over the equi-dimensional decomposition
\eqref{eq:equiDec}, where \( k \) takes the value \( k_i \) in \( \Omega_i^\varepsilon \) for \( i \in \{1, 2, f\} \),
and is completed with Neumann and Dirichlet boundary conditions. The problem can be reformulated as follows:
\begin{alignat}{2}
\nabla \cdot \mathbf{u}_i^\varepsilon & = 0  & \quad & \text{in } \Omega_i^\varepsilon,  \quad i = 1, 2, f, \label{eq:balancethin}\\
\mathbf{u}_i^\varepsilon & = -k_i \nabla p_i^\varepsilon & \quad & \text{in } \Omega_i^\varepsilon,  \quad i = 1, 2, f, 
\label{eq:closurethin} \\
p_i^\varepsilon & = p_f^\varepsilon & \quad & \text{on } \gamma_i,  \quad i = 1, 2, \label{eq:continuitypeq} \\
\mathbf{u}_i^\varepsilon \cdot \mathbf{n}_i^\varepsilon & = -\mathbf{u}_f^\varepsilon \cdot \mathbf{n}_f^\varepsilon & \quad & \text{on } \gamma_i,  \quad i = 1, 2, \label{eq:continuityueq}\\
-\mathbf{u} \cdot \mathbf{n} & = h &  \quad &\text{on } \Gamma^\mathbf{u},\\
p & = g & \quad & \text{on } \Gamma^p. \label{eq:dirCond3}
\end{alignat}
Conditions~\eqref{eq:continuitypeq} and~\eqref{eq:continuityueq} enforce the continuity of \( p \) and the normal flux across the interface \( \gamma_i \).

We derive interface conditions by integrating the balance law and the closure relation 
in the fracture along the segments \( [\mathsf{S}_1, \mathsf{S}_2] \), where 
\( \mathsf{S}_1 = \mathsf{S} - \frac{\varepsilon(\mathsf{S})}{2} \boldsymbol{\eta} \) 
and \( \mathsf{S}_2 = \mathsf{S} + \frac{\varepsilon(\mathsf{S})}{2} \boldsymbol{\eta} \). 
For simplicity, we neglect the dependence of \( \boldsymbol{\eta} \) on \( \mathsf{S} \) and of 
\( \mathbf{n}_i^\varepsilon \) on \( \mathsf{S}_i \) (for \( i \in \{1, 2\} \), while explicitly denoting 
the dependence of \( \mathbf{n}_f^\varepsilon \) on \( \mathsf{S}_i \).

\subsection{Statement}
The following assumptions are made in the derivation:
\begin{itemize}
    \item[A1] Functions $$u_{f,\eta}( \mathsf{S}+ r \frac{\varepsilon}{2}\boldsymbol{\eta}) \text{ and } \nabla_\tau \cdot k_f \nabla_\tau p_f^\varepsilon ( \mathsf{S}+ r \frac{\varepsilon}{2}\boldsymbol{\eta})$$ are linear in $r \in (-1,1)$;
        \item[A2] \( \mathbf{n}_f^\varepsilon(\mathsf{S}_1) =-\mathbf{n}_1^\varepsilon \approx -\boldsymbol{\eta} \) and 
              \( \mathbf{n}_f^\varepsilon(\mathsf{S}_2) =-\mathbf{n}_2^\varepsilon \approx \boldsymbol{\eta} \);
    \item[A3] \( \nabla_\tau \varepsilon \cdot \nabla_\tau \jump{p} \approx 0 \).
\end{itemize}

Under these assumptions, the problem with thin inclusions, given by 
Equations~\eqref{eq:balancethin}--\eqref{eq:dirCond3}, can be reformulated 
as an interface problem with Wentzell interface conditions for the flux jump 
and Robin interface conditions for the flux average. Specifically, the coefficients 
in the interface conditions~\eqref{eq:wentjump} and~\eqref{eq:wentmean} are:
\[
a^j = \varepsilon k_f, \quad r^m = \frac{k_f}{\varepsilon}, \quad r^j = h^j = a^m = g^m = 0.
\]

We note that A2 and A3 are satisfied for inclusions with constant aperture. 
A similar assumption to A3 is also commonly made in standard hybrid-dimensional 
models for fractured media~\cite{boon2018robust}.

\subsection{Steps in the Derivation}
To derive the interface conditions, we proceed as follows:

\noindent
\textbf{1.~Reformulation of the problem in the inclusion}\\
Exploiting the decomposition into normal and tangential operators of the divergence~\eqref{eq:divEq} 
and using~\eqref{eq:normalDiv}, 
we see that the balance equation can be written as:
\begin{equation}
\label{step0_1}
\frac{\partial u_{f,\eta}^\varepsilon}{\partial \eta} + \nabla_\tau \cdot \mathbf{u}_{f,\tau}^\varepsilon = 0.
\end{equation}
Instead, exploiting the decomposition for the closure relation, we have:
\begin{alignat}{2}
\mathbf{u}_{f,\tau}^\varepsilon & = - k_f \nabla_\tau p_f^\varepsilon, \label{step0_3} \\
u_{f,\eta}^\varepsilon & = -k_f \frac{\partial p_f^\varepsilon}{\partial \eta}. \label{step0_2}
\end{alignat}
Replacing Equation~\eqref{step0_3} into \eqref{step0_1}, we obtain:
\begin{alignat}{2}
\frac{\partial u_{f,\eta}^\varepsilon}{\partial \eta} & = \nabla_\tau \cdot k_f \nabla_\tau p_f^\varepsilon, \label{step1_1} \\
u_{f,\eta}^\varepsilon & = -k_f \frac{\partial p_f^\varepsilon}{\partial \eta}. \label{step1_2}
\end{alignat}

\noindent
\textbf{2.~Integration of the equations across the inclusion}\\
We integrate Equations~\eqref{step1_1} and \eqref{step1_2} along the segment \((\mathsf{S}_1, \mathsf{S}_2)\). While the left-hand side of Equation~\eqref{step1_1} and the right-hand side of Equation~\eqref{step1_2} can be explicitly evaluated, the other two terms can be integrated using assumption A1, yielding:
\begin{alignat}{2}
u_{f,\eta}^\varepsilon(\mathsf{S}_2) - u_{f,\eta}^\varepsilon(\mathsf{S}_1) & = \frac{\varepsilon}{2} \nabla_\tau \cdot k_f \nabla_\tau \Big( p_f^\varepsilon (\mathsf{S}_1) + p_f^\varepsilon (\mathsf{S}_2) \Big), \label{step4_1} \\
\frac{\varepsilon}{2} \Big( u_{f,\eta}^\varepsilon(\mathsf{S}_1) + u_{f,\eta}^\varepsilon(\mathsf{S}_2) \Big) & = -k_f \Big( p_f^\varepsilon(\mathsf{S}_2) - p_f^\varepsilon(\mathsf{S}_1) \Big). \label{step4_2}
\end{alignat}

\noindent
We observe that assumption A1 is equivalent to employing a trapezoidal quadrature rule for the integration of the two terms.


\noindent
\textbf{3.~Imposing continuity conditions}\\
Since \(u_{f,\eta}^\varepsilon(\mathsf{S}_i) = \mathbf{u}_f^\varepsilon(\mathsf{S}_i) \cdot \boldsymbol{\eta}\), 
exploiting assumption A2, the continuity conditions\eqref{eq:continuitypeq} and~\eqref{eq:continuityueq},
and A2 again, we have:
\[
\mathbf{u}_f^\varepsilon(\mathsf{S}_1) \cdot \boldsymbol{\eta} 
\approx -\mathbf{u}_f^\varepsilon(\mathsf{S}_1) \cdot \mathbf{n}_f^\varepsilon(\mathsf{S}_1)
= \mathbf{u}_1^\varepsilon(\mathsf{S}_1) \cdot \mathbf{n}_1^\varepsilon 
\approx \mathbf{u}_1^\varepsilon(\mathsf{S}_1) \cdot \boldsymbol{\eta},
\]
and
\[
\mathbf{u}_f^\varepsilon(\mathsf{S}_2) \cdot \boldsymbol{\eta} 
\approx \mathbf{u}_f^\varepsilon(\mathsf{S}_2) \cdot \mathbf{n}_f^\varepsilon(\mathsf{S}_2)
= -\mathbf{u}_2^\varepsilon(\mathsf{S}_2) \cdot \mathbf{n}_2^\varepsilon 
\approx \mathbf{u}_2^\varepsilon(\mathsf{S}_2) \cdot \boldsymbol{\eta}.
\]
Hence, Equations~\eqref{step4_1} and \eqref{step4_2} become:
\begin{alignat}{2}
\Big( \mathbf{u}_2^\varepsilon(\mathsf{S}_2) - \mathbf{u}_1^\varepsilon(\mathsf{S}_1) \Big) \cdot \boldsymbol{\eta} & = \frac{\varepsilon}{2} \nabla_\tau \cdot k_f \nabla_\tau \Big( p_1^\varepsilon(\mathsf{S}_1) + p_2^\varepsilon(\mathsf{S}_2) \Big), \label{step6_1} \\
\frac{\varepsilon}{2} \Big( \mathbf{u}_1^\varepsilon(\mathsf{S}_1) + \mathbf{u}_2^\varepsilon(\mathsf{S}_2) \Big) \cdot \boldsymbol{\eta} & = -k_f \Big( p_2^\varepsilon(\mathsf{S}_2) - p_1^\varepsilon(\mathsf{S}_1) \Big). \label{step6_2}
\end{alignat}

\noindent
\textbf{4.~Limit to interface}\\
We take the limits \(\mathsf{S}_i \to \mathsf{S}\) for \(i=1,2\) in Equations~\eqref{step6_1} and \eqref{step6_2}, and, employing the notation of jump and average, we get:
\begin{alignat}{2}
\jjump{\mathbf{u}} & = \varepsilon \nabla_\tau \cdot k_f \nabla_\tau \mean{p}, \label{step8_1} \\
\mmean{\mathbf{u}} & = -\frac{k_f}{\varepsilon} \jump{p}, \label{step8_2}
\end{alignat}
where we removed the superscript \(^\varepsilon\), as the functions are now defined on the subdomains \(\Omega_1\) and \(\Omega_2\).

Finally, exploiting assumption A3, the interface conditions~\eqref{step8_1} and \eqref{step8_2} can be written as:
\begin{alignat}{2}
-\jjump{\mathbf{u}} & = -\nabla_\tau \cdot (\varepsilon k_f \nabla_\tau \mean{p}), \label{step8_1_final} \\
-\mmean{\mathbf{u}} & = \frac{k_f}{\varepsilon} \jump{p}. \label{step8_2_final}
\end{alignat}

We do not explicitly report the weak formulation associated with these interface conditions,  
as in the next section, we will derive it for a generic network of interfaces.


\section{The Case of Multiple and Intersecting Heterogeneities}

We now consider a family of \(n_F\) thin heterogeneities, defined as in Subsection~\ref{sec:dec} by interfaces \(\gamma^j\) and apertures \(\varepsilon^j\), with \(j \in \{1, 2, \ldots, n_F\}\). These interfaces decompose \(\Omega\) into \(n_S\) connected subdomains \(\Omega_i\). Examples are illustrated in Figures~\ref{fig:domain1} and~\ref{fig:domain2}. Each subdomain is characterized by a permeability \(k_i\), while the thin heterogeneities are characterized by permeabilities \(k_f^j\). 

Exploiting the derivation of the interface conditions in Section~\ref{sec:derivation} to each heterogeneity, we obtain an interface problem of the form:
\begin{alignat}{2}
-\nabla \cdot \mathbf{u}_i &= 0 & \quad & \text{in } \Omega_i, \quad \text{for } i = 1, 2, \ldots, n_S, \label{eq:balanceIntM} \\
\mathbf{u}_i &= -k_i \nabla p_i & \quad & \text{in } \Omega_i, \quad \text{for } i = 1, 2, \ldots, n_S, \label{eq:closureIntM}\\
-\jjump{\mathbf{u}} & = -\nabla_\tau \cdot \varepsilon^j k_f^j \nabla_\tau \mean{p} & \quad & \text{on } \gamma^j, \quad \text{for } j = 1, 2, \ldots, n_F, \label{eq:wentzellIntM}\\
-\mmean{\mathbf{u}} & = \frac{k_f^j}{\varepsilon^j} \jump{p} & \quad & \text{on } \gamma^j, \quad \text{for } j = 1, 2, \ldots, n_F, \label{eq:RobinIntM}\\
-\mathbf{u} \cdot \mathbf{n} &= h & \quad & \text{on } \Gamma^{\mathbf{u}}, \label{eq:NeuIntM}\\
p &= g & \quad & \text{on } \Gamma^p. \label{eq:DirIntM}
\end{alignat}

By extending the definition of the function space \(V\) and the function set \(U\) from Eq.~\eqref{eq:defV} to multiple subdomains, we derive the weak formulation in the usual manner. From the left-hand side of Eq.~\eqref{eq:balanceIntM}, we have:
\[
\begin{aligned}
\sum_{i=1}^{n_S} \int_{\Omega_i} \nabla \cdot \mathbf{u}_i \, q_i \, \mathrm{d}\Omega &= \sum_{i=1}^{n_S} \left( - \int_{\Omega_i} \mathbf{u}_i \cdot \nabla q_i \, \mathrm{d}\Omega + \int_{\partial \Omega_i} \mathbf{u}_i \cdot \mathbf{n}_i \, q_i \, \mathrm{d}\Gamma \right) \\
&= \sum_{i=1}^{n_S} \left( - \int_{\Omega_i} \mathbf{u}_i \cdot \nabla q_i \, \mathrm{d}\Omega \right) + \int_{\Gamma} \mathbf{u}_i \cdot \mathbf{n}_i \, q_i \, \mathrm{d}\Gamma \\
&\quad + \sum_{j=1}^{n_F} \left( - \int_{\gamma^j} \jjump{\mathbf{u}} \mean{q} \, \mathrm{d}\Gamma - \int_{\gamma^j} \mmean{\mathbf{u}} \jump{q} \, \mathrm{d}\Gamma \right),
\end{aligned}
\]
where \(q\) is a test function in \(V\).

In the last term above, substituting Darcy's law~\eqref{eq:closureIntM}, the interface conditions~\eqref{eq:wentzellIntM} and~\eqref{eq:RobinIntM}, and the Neumann boundary conditions~\eqref{eq:NeuIntM}, leads to:
\[
\begin{aligned}
\sum_{i=1}^{n_S} \int_{\Omega_i} k_i \nabla p_i \cdot \nabla q_i \, \mathrm{d}\Omega 
&+ \sum_{j=1}^{n_F} \Bigg( - \int_{\gamma^j} \nabla_\tau \cdot \big(\varepsilon^j k_f^j \nabla_\tau \mean{p}\big) \mean{q} \, \mathrm{d}\Gamma \\
&\quad + \int_{\gamma^j} \frac{k_f^j}{\varepsilon^j} \jump{p} \jump{q} \, \mathrm{d}\Gamma \Bigg) = \int_{\Gamma^{\mathbf u}} h \, q \, \mathrm{d}\Gamma.
\end{aligned}
\]
Finally, applying integration by parts to the diffusion terms on \(\gamma^j\) and exploiting the hypothesis~\eqref{eq:hp}, we obtain the weak formulation:

\vspace{0.5em}
\noindent
Find \(p \in U\) such that:
\[
\begin{aligned}
\sum_{i=1}^{n_S} \int_{\Omega_i} k_i \nabla p_i \cdot \nabla q_i \, \mathrm{d}\Omega 
&+ \sum_{j=1}^{n_F} \Bigg(  \int_{\gamma^j} \varepsilon^j k_f^j \nabla_\tau \mean{p} \cdot \nabla_\tau \mean{q} \, \mathrm{d}\Gamma \\
&\quad + \int_{\gamma^j} \frac{k_f^j}{\varepsilon^j} \jump{p} \jump{q} \, \mathrm{d}\Gamma \Bigg) = \int_{\Gamma^{\mathbf u}} h \, q \, \mathrm{d}\Gamma \quad \text{for all } q \in V.
\end{aligned}
\]

This resulting formulation is elegant and simple, involving two types of interface integrals: one arising from the Wentzell interface conditions for the flux jump and the other from the Robin interface condition for the flux average. It can easily accommodate any fracture network, where fractures may have different conductivities and apertures.

\subsection{Computing the Solution at Interfaces}

In this modeling approach, the equations within the fractures are incorporated  
through interface conditions, with the solution computed only in the background domain.  
However, in many applications, it is of interest to evaluate the pressure along the fractures.  
This can be achieved using, for example, any of the relations presented in~\cite{martin2005modeling}.  

The method we propose is to use the average of \(p\) at the interface as the pressure in the fracture:  
\begin{equation}
\label{eq:pf}
p_f = \mean{p}.
\end{equation}

Using this relation, we observe that the proposed method is closely related to  
averaging-based hybrid-dimensional models.  
In fact, if we adopt the relation~\eqref{eq:pf} as a definition  
and substitute it into~\eqref{eq:wentzellIntM},  
we obtain an expression that is remarkably similar to that of averaging-based hybrid-dimensional models.  

Therefore, the presented approach can be interpreted as a reformulation of such models,  
where \(p_f\) is formally eliminated in a manner analogous to static condensation.

\subsection{Finite Element Formulation}

Moreover, the formulation is readily applicable in finite element methods. Specifically, we can proceed in the standard way by introducing a mesh of triangles or quadrilaterals in 2D, or tetrahedra or hexahedra in 3D, for each subdomain \(\Omega_i\), ensuring that the boundaries of elements align with the interfaces \(\gamma^j\).
One approach to generate such a mesh is to create a global mesh for the whole domain \(\Omega\) that captures the geometry of \(\gamma^j\) and then splitting it into the different subdomains. The boundaries of the mesh elements at the interfaces are duplicated, creating a structure analogous to domain decomposition methods.

\section{Numerical validation}

We validate the proposed model for interface conditions through four test cases of fluid flow in fractured media. The first three test cases, which vary in complexity, were introduced in previous studies as benchmark problems to compare different discretization methods. The first benchmark is the \emph{regular fracture network (2D)} presented in~\cite{flemisch2018benchmarks}. The second and third benchmarks are the \emph{Regular Fracture Network (3D)} and the \emph{single fracture}, as described in~\cite{BenchPaper}. Since these three test cases involve fractures with a constant aperture, we introduce a fourth test case: a two-dimensional problem featuring a fracture with an elliptical shape.
Validation is performed by comparing pressure profiles along selected segments, as proposed in the two benchmark articles.

Regarding the discretization method, we employ a standard finite element method using linear or bilinear elements, depending on the mesh used. As previously described, the employed meshes resolve the fractures by cutting along the elements. This approach ensures that the resulting meshes after the split remain conforming on both sides of each interface.

The test cases were implemented in {\tt parroth}, an application developed within the FE framework {\tt MOOSE}, available at the repository {\tt https://github.com/favinom/parroth}. The {\tt MOOSE} framework natively provides functionalities for mesh splitting through the {\tt BreakMeshByBlockGenerator} function and supports handling interface problems. All linear systems resulting from the FE discretization of these problems were solved using the parallel direct solver {\tt MUMPS}.

\subsection{Regular fracture network (2D)}

\begin{figure}
\includegraphics[scale=0.4]{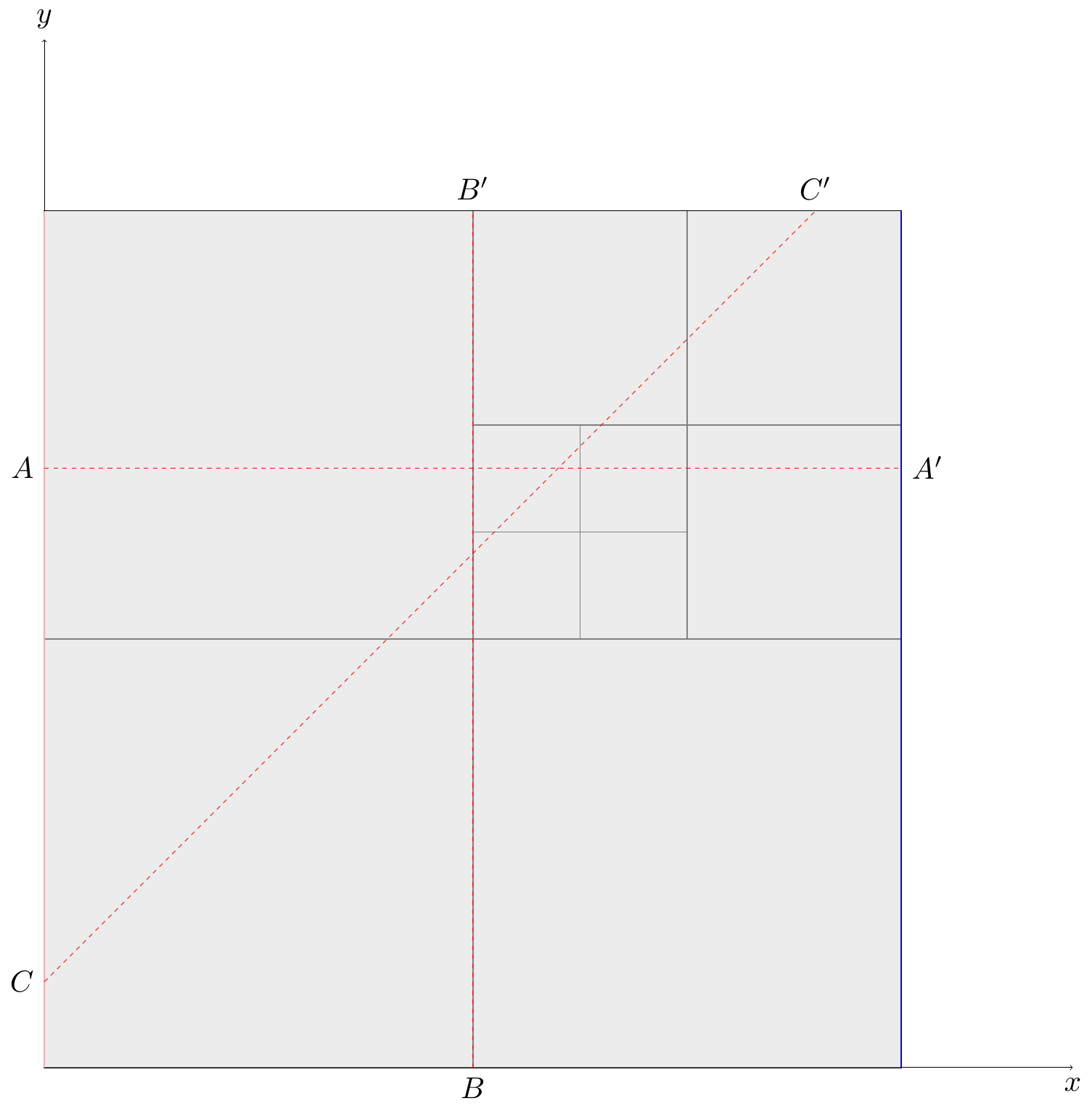}
\caption{Domain and fracture network for the test case \emph{Regular Fracture Network (2D)}. The left side, in pink, is where flow boundary conditions are imposed, while the right side, in blue, is where pressure boundary conditions are applied.
Pressure profiles for the conductive scenario are evaluated along the segments $AA'$ and $BB'$, while a pressure profile for the blocking scenario is evaluated along the segment $CC'$.}

\label{fig:domain1}
\end{figure}

The domain of this test case is the square $(0,1)^2$,
where six fractures divide it into ten subdomains, as shown in Figure~\ref{fig:domain1}, which also highlights the segments along which the pressure profiles are evaluated. The permeability in each subdomain of the background is set to $k_i = 1$, with $i\in\{ 1,2,\ldots 10\}$, and all fractures have an aperture of $\epsilon=10^{-4}$. Two different scenarios are considered for the fracture permeability: a conductive scenario, where all fractures feature a permeability $k_f^j=10^4$, and a blocking scenario, where all fractures feature a permeability $k_f^j=10^{-4}$.
Pressure boundary conditions ($p=1$) are imposed on the right side of the domain at $x=1$, reported in light blue in Figure~\ref{fig:domain1}, while flux boundary conditions ($h=1$) are imposed on the left side of the domain at $x=0$, reported in pink in Figure~\ref{fig:domain1}.
On the top and bottom sides, homogeneous flux boundary conditions are imposed.

For the discretization, a uniform, regular mesh with 32 square elements per side has been employed. Hence, the total number of mesh vertices after duplicating the ones at the interfaces is 1210, which is in line with the other methods employed in the benchmark paper. In fact, only the BOX method~\cite{helmig1997multiphase} presents a smaller number of nodes than ours, while the other methods presented a larger number.

\subsubsection{Conductive scenario}

\begin{figure}
\includegraphics[trim={19cm 5cm 12cm 5cm},clip,scale=0.2]{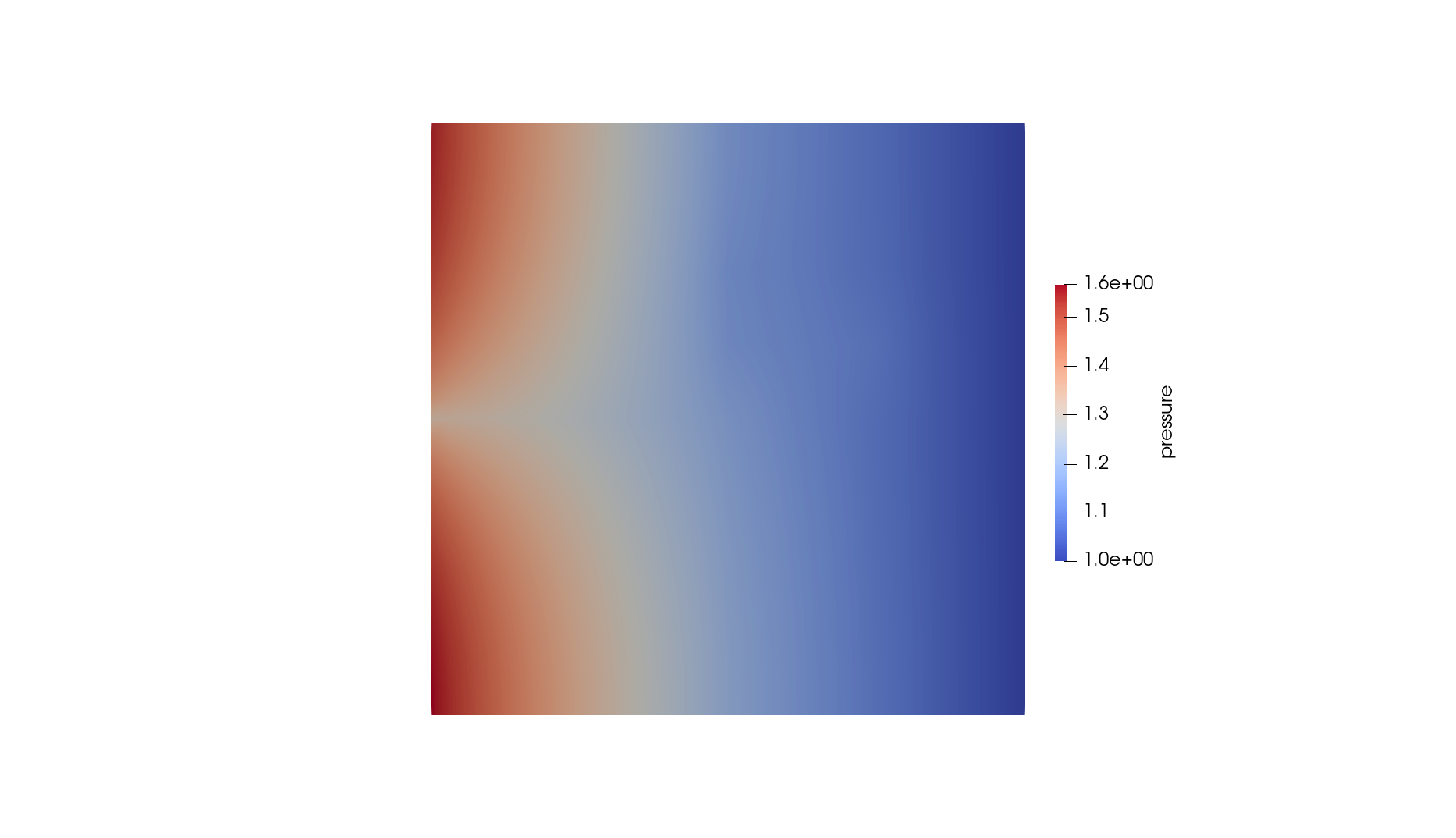}
\caption{Pressure distribution for the conductive scenario ($k_f=10^4$) of the test case \emph{Regular Fracture Network (2D)}.}
\label{fig:pressureCond2D}
\end{figure}

The pressure distribution obtained for the conductive scenario is reported in Figure~\ref{fig:pressureCond2D}, which is graphically identical to the one reported in the reference article.
For this scenario, the pressure distribution along two segments is reported: $AA'$ at $y=0.7$ and $BB'$.
This latter is actually a profile along the vertical fracture at $x=0.5$.

%

\begin{figure}
\includegraphics[scale=0.4]{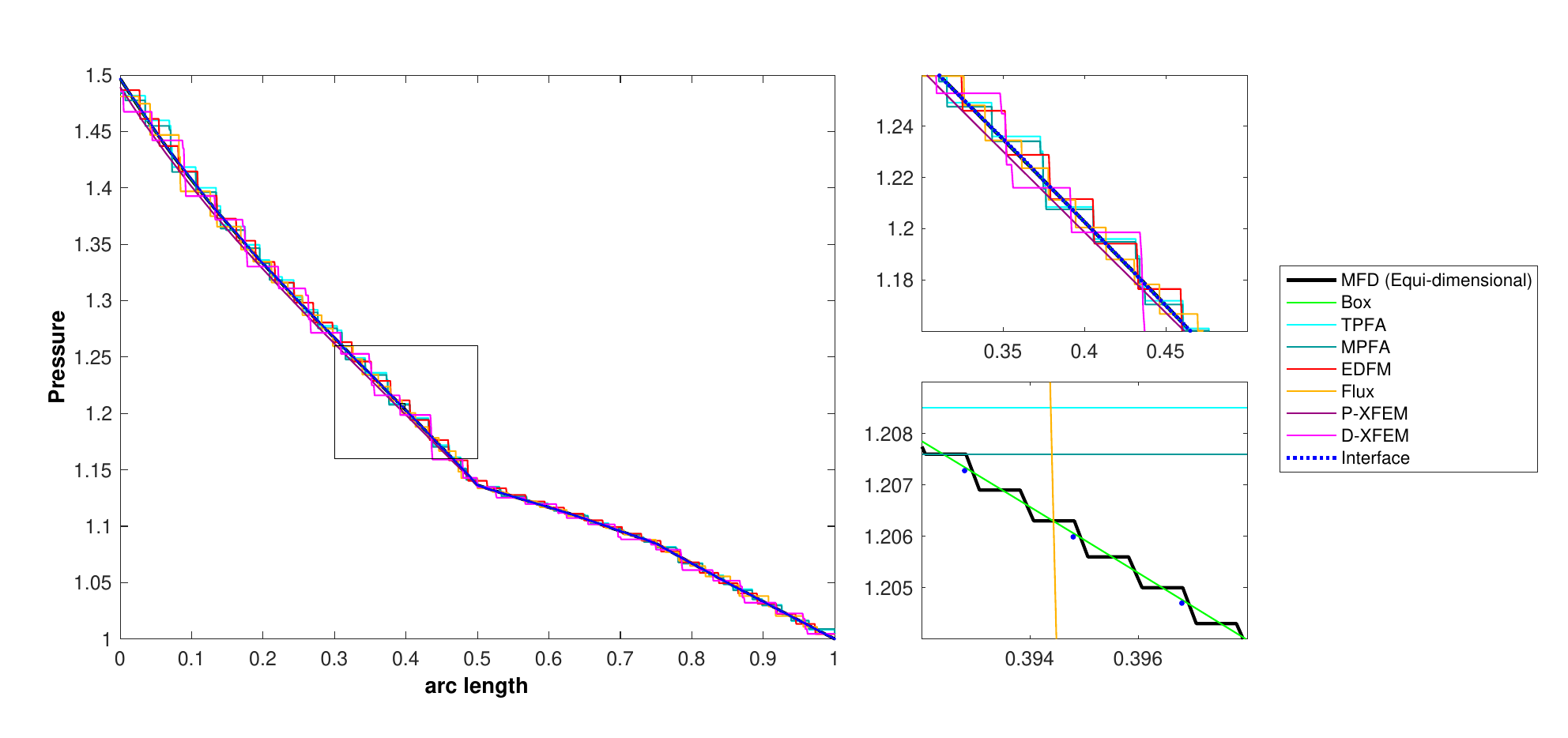}
\caption{Comparison of pressure profile along the segment $AA'$ at $y=0.7$.
In the main plot along the whole segment, the reference solution is hidden by the one obtained with the BOX method and by our solution obtained with the model based on interface conditions, reported in dotted blue line.
}
\label{fig:profileCondY}
\end{figure}

A comparison between the results obtained with a FE discretization
of the proposed interface model and the methods presented in~\cite{}
is reported in Figure~\ref{fig:profileCondY}.
Our results, like those presented in the article, are close to the reference solution, which has been computed using a mimetic finite differences (MFD) method applied to the equi-dimensional model employing a extremely fine grid with more than two millions degrees of freedom.
No notable differences are observed between methods based on hybrid-dimensional models, whether they are continuity-based or use immersed approaches.
In the main plot, our results are superimposed with those obtained using the BOX method, which employs an immersed approach, and both models are superimposed with the reference solution.
In the zoomed-in plots, minor differences can be observed between the reference solution and our model. Upon closer inspection of the pressure profile, these differences appear to result from interpolation characteristics used to report the error in the reference solution.

\begin{figure}
\includegraphics[scale=0.4]{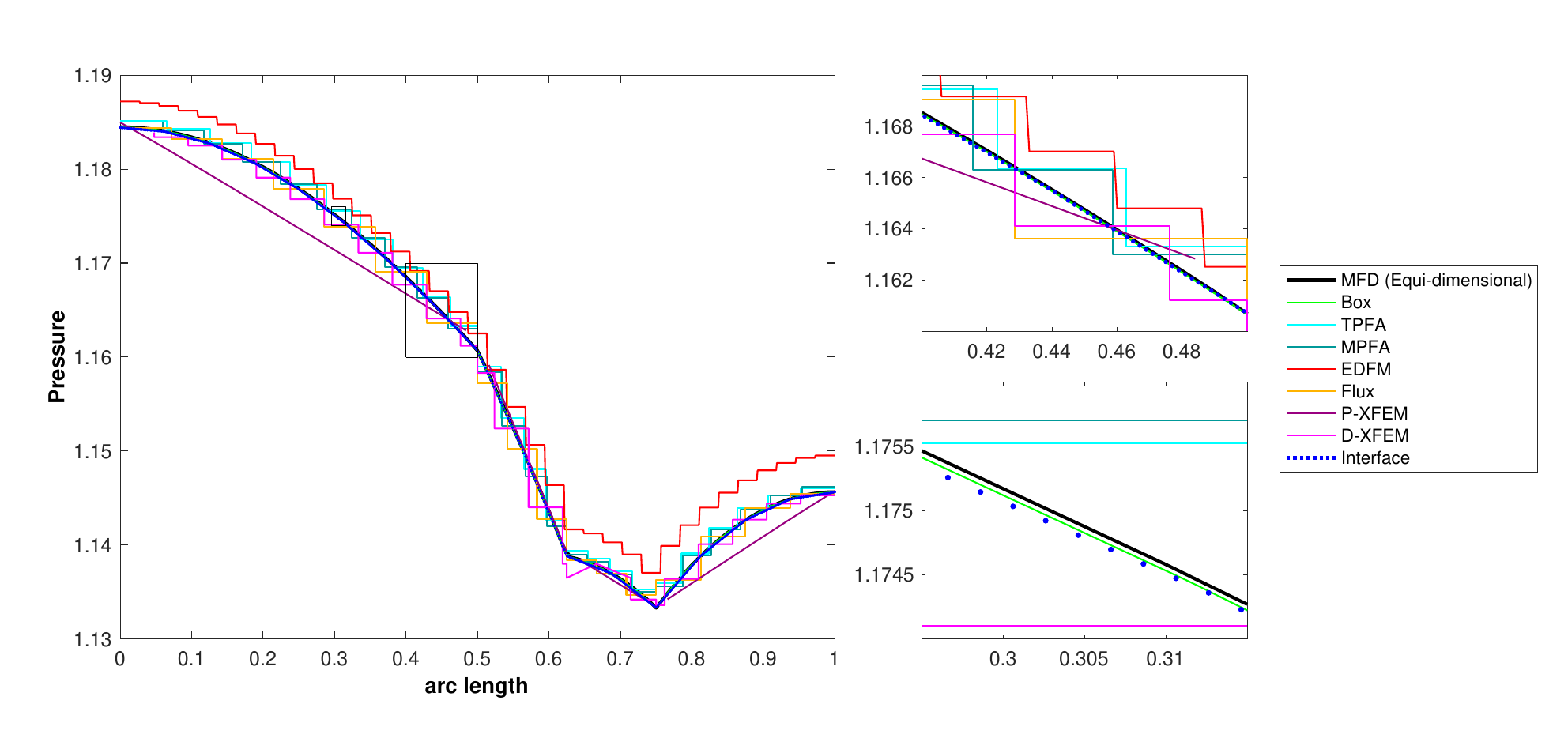}
\caption{Comparison of pressure profile along the segment $BB'$ at $x=0.5$.
In the main plot along the whole segment, the reference solution is hidden by the one obtained with the BOX method and by our solution obtained with the model based on interface conditions, reported in dotted blue line.
}
\label{fig:profileCondX}
\end{figure}

In Figure~\ref{fig:profileCondX}, the same comparison is performed along the vertical fracture at $x=0.5$
so that we computed it exploiting the formula~\eqref{eq:pf}.
Unlike the pressure profiles obtained along the segment $AA'$,
the various models and methods show some significant deviations from the reference solution.
However, the results obtained with a finite element discretization of our interface model
are perfectly superimposed with both the reference solution and the BOX method,
which is based on a continuous pressure model.

\subsubsection{Blocking scenario}

\begin{figure}
\includegraphics[trim={19cm 5cm 12cm 5cm},clip,scale=0.2]{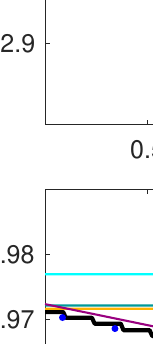}
\caption{Pressure distribution for the blocking scenario ($k_f=10^{-4}$) of the test case \emph{Regular Fracture Network (2D)}.}
\label{fig:pressureBlock2D}
\end{figure}

\begin{figure}
\includegraphics[scale=0.4]{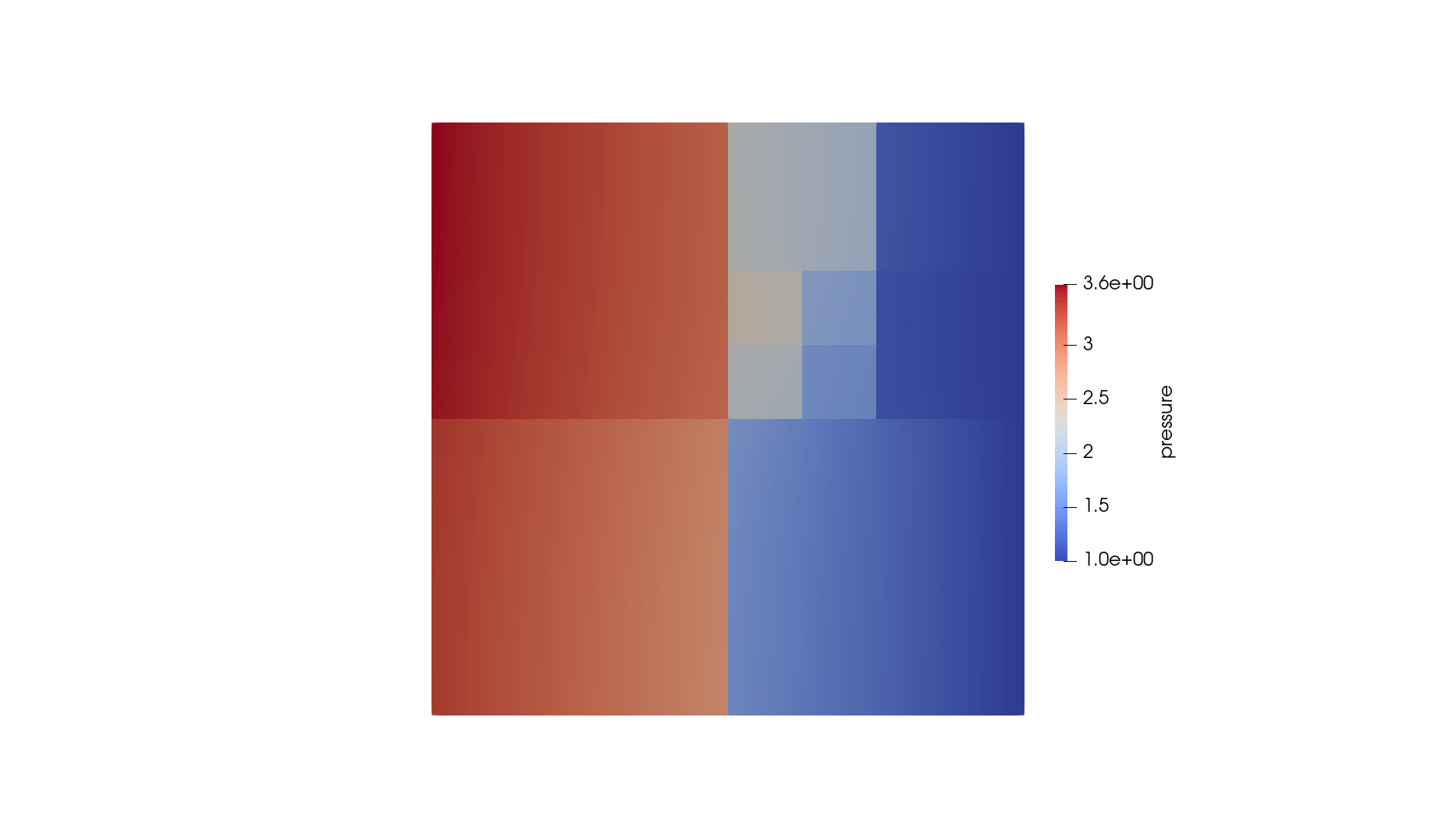}
\caption{Comparison of pressure profile along the segment $CC'$.
In the main plot along the whole segment, the reference solution is hidden by the one obtained with the P-XFEM method and by our solution obtained with the model based on interface conditions, reported in dotted blue line.
}
\label{fig:profileBlockX}
\end{figure}

The pressure distribution for the blocking scenario is shown in Figure~\ref{fig:pressureBlock2D}. 
In this scenario, the fractures exhibit a permeability lower than that of the surrounding medium, resulting in pressure discontinuities at the fracture interfaces. Traditional models that assume continuous pressure generally fail to capture this discontinuity, as illustrated in Figure~\ref{fig:profileBlockX}. In contrast, averaging-based hybrid-dimensional models are generally more accurate.
Our proposed model, which is based on interface conditions, successfully reproduces the expected pressure profile. The results align closely with the reference solution and the P-XFEM (primal XFEM) method, indicating accurate handling of the pressure discontinuity.

As a conclusion from the numerical simulation of this first test case, we observe that our model, based on interface conditions, effectively handles the intersection of multiple fractures without requiring specific fracture equations or meshes for the fractures and their intersections.

\subsection{Regular fracture network (3D)}

\subsection{Single fracture (3D)}

\subsection{Elliptical Fracture (2D)}

In this final numerical test case,  
we evaluate the ability of our model to handle fractures with varying aperture.  
The domain is \(\Omega = (0,1)^2\), containing a single vertical elliptical fracture  
centered at \((0.5,0.5)\), with a minor axis of size \(10^{-4}\) and a major axis of size \(1 + 10^{-4}\).  
The permeabilities are set to \(k_1 = k_2 = 1\) and \(k_f = 10^{-4}\).  
We compare the pressure distribution obtained from the equi-dimensional model  
and the reduced model based on interface conditions.  

For the equi-dimensional model,  
the mesh contains 1002 elements in the \(x\)-direction and 1000 elements in the \(y\)-direction.  
The elements along the \(x\)-axis are adjusted to conform to the interface boundary,  
ensuring two elements across the fracture aperture.  
For the reduced model with interface conditions,  
we use a uniform mesh with 1000 elements in both directions,  
splitting the domain along \(x = 0.5\).  

\begin{figure}
\includegraphics[scale=0.2]{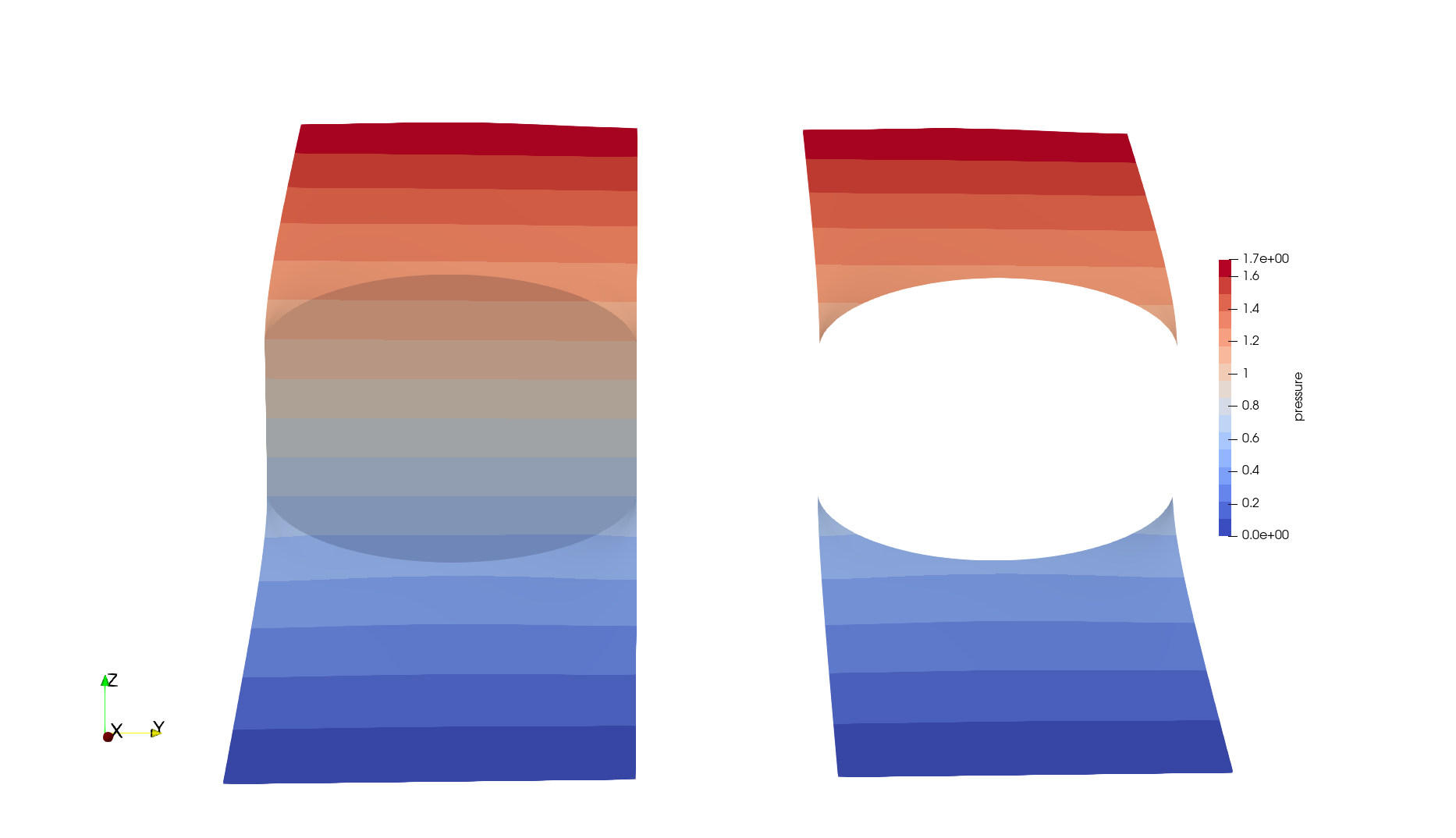}
\caption{Comparison between the equi-dimensional model and the interface-condition-based model for a vertical elliptical fracture.}
\label{fig:ellipseComparison}
\end{figure}

Figure~\ref{fig:ellipseComparison} presents the comparison between the equi-dimensional model  
and the model based on interface conditions.  
The results are macroscopically similar, validating the assumption A3 of small aperture variations.  
For the equi-dimensional model, the pressure drop along the \(x\)-axis  
is influenced by the fracture aperture, and the model based on interface conditions  
successfully reproduces this behavior through a discontinuity depending on the point \(\mathsf{S} \in \Gamma\) along the interface.  
This is further confirmed by additional analyses comparing profiles along segments (not shown here).  
These findings demonstrate that, under the assumption A3, the proposed model effectively manages fractures  
with varying aperture directly.

\section{Conclusions}

We presented a reduced modeling technique in which  
a suitable interface problem is used to approximate  
a heterogeneous diffusion problem,  
written in the form of divergence of a flux,  
featuring thin heterogeneities,  
such as the equations describing the fluid flow in fractured porous media.  
By performing a formal integration of the equations across the aperture of the inclusions,  
conditions for the flux jump and average at the interface naturally emerge.  
In particular, we showed that the jump of the flux  
is equal to a tangential diffusion operator applied to the average of the solution across the interface,  
resulting in a sort of Wentzell-type interface condition,  
while for the average of the flux, the conditions are of Robin type.  

The resulting interface conditions can be straightforwardly incorporated into the weak  
formulation of the corresponding interface problem,  
resulting in an approach that is exceptionally elegant and can be readily employed for finite element discretizations.  
In particular, unlike other reduced modeling techniques based on hybrid-dimensional representations of thin inclusions,  
the proposed technique eliminates the need to explicitly account for the equations within the inclusions.  

We demonstrated the accuracy and robustness of this modeling technique through extensive numerical validation.  
We showed that it can effectively handle a wide range of material properties,  
naturally accommodate networks of inclusions without the need for an explicit mesh,  
and manage fractures with varying apertures.  

While we focused on heterogeneous diffusion problems with a scalar diffusion coefficient  
and inclusions that decompose the domain of the problem,  
this modeling technique can be adapted  
to handle tensorial material properties and fully embedded inclusions,  
and it is adaptable to general problems in divergence form.
 
The flexibility of this technique opens up new avenues for future research.
At the modeling level, its application to more complex and multiphysics problems,
such as elasticity and poroelasticity, holds significant promise but would require thorough validation.
On the analytical side, the regularity of the solution remains a critical question.
While we expect that regularity results for Wentzell boundary conditions
can be extended to cases involving inclusions that separate the domain into connected subsets,
further analysis is needed to establish the regularity of solutions involving fully embedded inclusions.
From a numerical perspective, this modeling technique is particularly well-suited for advanced discretization methods,
such as extended finite element methods (X-FEM) or cut finite element methods (CutFEM).
The integration of these methods with the proposed modeling technique
would enable the seamless treatment of networks of inclusions,
eliminating the need for meshes that conform to the geometry of the inclusions.

\bibliographystyle{elsarticle-num}
\bibliography{references}

\end{document}